%% file: p-adic-height-computation.tex
\documentclass{tran-l}[12pt]
\usepackage[all]{xy}
\input macros.tex

\input defs.tex
\begin{document}
\title{Computing local $p$-adic height pairings on hyperelliptic curves}

\author{Jennifer S. Balakrishnan}
\address{
Department of Mathematics\\
Massachusetts Institute of Technology\\
77 Massachusetts Avenue\\
Cambridge, MA 02139
}

\author{Amnon Besser}
\address{
Department of Mathematics\\
Ben-Gurion University of the Negev\\
P.O.B. 653\\
Be'er-Sheva 84105\\
Israel
}

\maketitle

\input intro.tex 
\input coleman-gross.tex 
\input hyperelliptic.tex 
\input coleman-integration.tex 

\input algorithm.tex 
\input implementation.tex 

\bibliographystyle{amsalpha}
\bibliography{biblio} 

\end{document}

%% file: macros.tex
\usepackage{amssymb}

\renewcommand{\div}{\mbox{\rm div}}
\DeclareMathOperator{\chr}{char}

\DeclareMathOperator{\Pic}{Pic}
\DeclareMathOperator{\Lie}{Lie} 
\DeclareMathOperator{\Frob}{Frob}

\DeclareMathOperator{\Res}{Res}

\newcommand{\inv}{^{-1}}

\newcommand{\ra}{\rightarrow}

\newcommand{\cO}{\mathcal{O}}

 \DeclareMathOperator{\Div}{Div}
 \DeclareMathOperator{\res}{res}

\newcommand{\lra}{\longrightarrow}

\newcommand{\Ga}{\mathbb{G}_a}

\newcommand{\comment}[1]{}
\newcommand{\Q}{\mathbb{Q}}

\newcommand{\C}{\mathbb{C}}

\newcommand{\Z}{\mathbb{Z}}

\renewcommand{\O}{\mathcal{O}}
\newcommand{\A}{\mathbb{A}}

 \DeclareMathOperator{\tr}{tr}

\theoremstyle{plain}
\newtheorem{theorem}{Theorem}[section]
\newtheorem{proposition}[theorem]{Proposition}
\newtheorem{corollary}[theorem]{Corollary}

\newtheorem{lemma}[theorem]{Lemma}

\theoremstyle{definition}
\newtheorem{definition}[theorem]{Definition}
\newtheorem{algorithm}[theorem]{Algorithm}

\theoremstyle{remark}

\newtheorem{remark}[theorem]{Remark}

\newtheorem{example}[theorem]{Example}

\numberwithin{equation}{section}
\numberwithin{figure}{section}
\numberwithin{table}{section}

\newcounter{listnum}


%% file: defs.tex
\renewcommand{\O}{\mathcal{O}}
\newcommand{\Rw}{^\textbf{w}}

\newcommand{\Rnw}{^\textbf{nw}}
\newcommand{\ndiv}{\nmid}

\newcommand{\dlog}{\operatorname{d\!\log}}

\newcommand{\lag}{\Psi}
\newcommand{\PP}{\mathbb{P}}
\newcommand{\dr}{\textup{dR}}
\newcommand{\hdr}{H_{\dr}}

\newcommand{\rev}{\textup{rev}}

\newcommand{\frev}{f^{\rev}}


%% file: intro.tex
\begin{abstract}
  We describe an algorithm to compute the local
  component at $p$ of the Coleman-Gross $p$-adic height pairing on
  divisors on hyperelliptic curves. As the height pairing is given in
  terms of a Coleman integral, we also provide new techniques to
  evaluate Coleman integrals of meromorphic differentials and present
  our algorithms as implemented in Sage.
\end{abstract}

\section{Introduction}
For an elliptic curve over $\Q$, the classical Birch and
Swinnerton-Dyer (BSD) conjecture predicts that a special value of its
$L$-function can be given in terms of certain arithmetic invariants of
the curve, one of which involves the canonical height pairing matrix
of a basis of rational points. The $p$-adic analogue \cite{mazur:padicBSD} of the
BSD conjecture makes a similar prediction, with the canonical height
pairing replaced by a $p$-adic one \cite{mazur-tate}. These
conjectures have natural generalizations to abelian varieties.

The $p$-adic height pairing was first defined by
Schneider~\cite{Schn82} for abelian
varieties and was extended to motives by
Nekov{\'a}\v{r}~\cite{Nek93}. For Jacobians of curves there is a third
definition, due to Coleman and Gross~\cite{coleman-gross}, relying on
Coleman's theory of $p$-adic
integration~\cite{coleman:dilogarithms,coleman:torsion,coleman-deshalit,besser:coleman}. This third definition of the height pairing
is known to be equivalent to the previous ones~\cite{Bes99a}.

For the purpose of numerically verifying $p$-adic BSD type conjectures, it is
important to have an effective algorithm for the computation of the
$p$-adic height pairing. By the work of
Kedlaya \cite{kedlaya:mw} and Mazur, Stein, and Tate
\cite{mazur-stein-tate}, we can easily compute $p$-adic heights on
elliptic curves. Our work deals with the next logical step, $p$-adic height pairings on Jacobians
of hyperelliptic curves.

The reason for treating Jacobians of hyperelliptic curves is that we
have available the Coleman-Gross definition of the height pairing,
which is much more concrete than previous definitions. The restriction
to hyperelliptic curves is made primarily so that we may apply
the recent algorithm~\cite{bbk} for the computation of Coleman
integrals on such curves, relying in turn on Kedlaya's work on the
computation of the matrix of Frobenius on hyperelliptic
curves~\cite{kedlaya:mw}. We note that generalizations of Kedlaya's work to other
types of curves can be applied to generalize the results to these
curves as well.

Coleman and Gross give a decomposition of the global $p$-adic height pairing as a sum of local
height pairings at each prime. The local heights away from the prime $p$ behave in much the
same way as local archimedean heights,  so the main interest lies in
the primes above $p$, where Coleman integration is used. It is this
last type of local height pairing which we aim to compute.

To fix ideas, consider a hyperelliptic curve $C$ over $\Q_p$ with
$p$ a prime of good reduction. Then, for $D_1,D_2 \in \Div^0(C)$ with
disjoint support, the Coleman-Gross $p$-adic height pairing at the
prime $p$ is given in terms of the Coleman integral
 $$h_p(D_1,D_2) = \int_{D_2} \omega_{D_1},$$ for
an appropriately constructed differential $\omega_{D_1}$ associated to
the divisor $D_1$. This last association is not straightforward and
relies again on Coleman integration.

We say next to nothing in this work about the computation of local
height pairings away from $p$. As mentioned before, this is not a
$p$-adic problem and is shared with the computation of archimedean
height pairings. 
In~\cite{besser:sage} we suggested a method for treating this
problem, which needed in particular some refined estimates of
Kausz~\cite{Kau99}. In the meantime, we have learned that this problem is
treated in the recent Ph.D. thesis of M\"uller~\cite{muller:thesis}.

The structure of the paper is as follows: in Section~\ref{sec:review},
we review the work of Coleman and Gross \cite{coleman-gross} that
constructs $p$-adic heights for curves, and in particular, defines the
local contribution at $p$ in terms of a Coleman integral. After a
brief review of hyperelliptic curves in Section~\ref{sec:hyper} we
give, in  Section~\ref{coleman}, an overview of Coleman integration, where we
discuss some known results for computing these integrals and describe
a new construction that allows us to compute a broader class of
integrals that the heights necessitate. We discuss the algorithm for
computing the local height pairings in Section~\ref{sec:localatp}.   In Section~\ref{notes}, we
discuss our implementation of the algorithm in Sage along with error
bounds on our results.  We follow this in Section~\ref{ex} with numerical examples illustrating our methods. We conclude in Section~\ref{future} by posing some questions arising from our work.

\subsubsection*{Acknowledgments.} JB would like to thank Kiran
Kedlaya for engaging in numerous
helpful discussions and Liang Xiao for much constructive conversation. The
implementation of this algorithm would not have been possible without
David Roe and Robert Bradshaw, who laid much of the groundwork in
Sage. Computations were done on \texttt{sage.math.washington.edu} (NSF
Grant No. DMS-0821725), thanks to the generosity of William Stein.
AB would like to thank the Institute for Advanced Study in Princeton
where some preliminary work on this project, leading to the very rough
sketch~\cite{besser:sage},
begun during the academic year 2006/2007, and the Isaac Newton Institute for
Mathematical Sciences in Cambridge where some further work on this
project was done in August 2009.
JB's work was made possible by the NDSEG and NSF Graduate
Fellowships. AB work was supported by  the Bell companies
Fellowship and the James
D. Wolfensohn fund while at the Institute for Advanced Study and is
presently supported by a grant from the Israel Science Foundation.

%% file: coleman-gross.tex
\section{The $p$-adic height pairing}
\label{sec:review}

In this section we review the definition of the Coleman-Gross height
pairing. As explained in the introduction, there are two required
ingredients for making this definition: the theory of Coleman
integration and a certain choice of a canonical form. These will be discussed
in detail in a later section.

Suppose $X/K$ is a curve defined over a number field $K$, with
good reduction at primes above $p$.
To define the height pairing
\begin{equation*}
  h: \Div^0(X) \times \Div^0(X) \to \Q_p\;,
\end{equation*}
one needs the following data:
\begin{itemize}
\item  A ``global log''- a continuous idele class character
 \begin{equation*}
   \ell: \mathbb{A}_K^\times/K^\times \to \Q_p\;.
 \end{equation*}
\item For each $v|p$ a choice of a subspace $W_v\in \hdr^1(X\otimes
  K_v/K_v)$ complementary to the space of holomorphic forms.
\end{itemize}
For the definition of the Coleman-Gross height we must insist that the
local characters $\ell_v$ induced by $\ell$, for $v|p$, are ramified
in the sense that they do not vanish on the units in $K_v$.

  From $\ell$ one deduces the following data:
\begin{itemize}
\item For any place $v\ndiv p$ we have $\ell_v(\O_{K_v}^\times)=0$ for
  continuity reasons, which implies that $\ell_v$ is completely
  determined by the number $\ell_v(\pi_v)$, where $\pi_v$ is any
  uniformizer in $K_v$.
\item For any place $v|p$ we can decompose
  \begin{equation*}
    \xymatrix{
      {\O_{K_v}^\times}  \ar[rr]^{\ell_v} \ar[dr]^{\log_v} & &   \Q_p\\
      & K_v\ar[ur]^{t_v}
    }
  \end{equation*}
  where $t_v$ is a $\Q_p$-linear map. Since we assume that $\ell_v$
  is ramified it is then possible to extend
  \begin{equation}\label{logbranch}
  \log_v : K_v^\times \to K_v
\end{equation}
 in such a way that the diagram remains commutative.
\end{itemize}

Let us now describe the height pairing $h(D_1,{D_2})$ for a pair of degree zero
divisors ${D_1}$ and ${D_2}$ with disjoint support. The height pairing is a sum of local
terms $$h({D_1},{D_2})= \sum_v h_v({D_1},{D_2})$$ over all finite places
$v$. The local terms depend only on the completion at $v$ of $K$. Thus,
let $k$ be a local field of characteristic 0, with valuation ring
$\cO$, uniformizer $\pi$ and let $F = \cO/\pi\cO$ be the residue
field, with order $q.$ Let $C$ denote the curve $X$ over the local field $k$. We shall assume that $C$ has a $k$-rational
point and that $C$ has good reduction at $\pi$.

Let $\chi: k^* \ra \Q_p$ be a continuous homomorphism, which is the
local component of $\ell$.

\begin{proposition}\label{intersect} If $\chr F \neq p$, and $D_1,D_2
  \in \Div^0(C)$ have disjoint support, then there is a unique function
  $\langle D_1,D_2 \rangle$ that is continuous, symmetric,
  bi-additive,taking values in $\Q_p$, and satisfying
  \begin{equation}
  \langle (f),  D_2\rangle = \chi(f(D_2))\label{eq:heighprinc}
\end{equation}
 for $f \in k(C)^*$.
\end{proposition}

\begin{proof}See \cite[Prop 1.2]{coleman-gross}.
In fact, one has~\cite[(1.3)]{coleman-gross}
\begin{equation}
  \label{locCG}
  h_v({D_1},{D_2}) = \ell_v(\pi_v) \cdot ({D_1},{D_2})\;.
\end{equation}
Here, $({D_1},{D_2})$ denotes intersection multiplicity on a regular model of
$C$ over $\O$ of extensions of ${D_1}$ and ${D_2}$ to this model. To make
this have the required properties one of these extensions has to have
zero intersection with all components of the special fibre.
\end{proof}

We now describe the local contribution at a place $v|p$.

\begin{definition}A meromorphic differential on $C$ over $k$ is said
  to be of the
  \emph{first kind} if it is holomorphic, of the \emph{second kind} if
  it has residue $0$ at every point, and of the \emph{third kind} if
  it has a simple
  pole with residue in $\Z$, respectively.
\end{definition}

Recall that the differentials of the second kind, modulo exact
differentials, i.e., differentials
of rational functions, form a finite dimensional $k$-vector space of
dimension
$2g$. It is canonically isomorphic
to the first algebraic de Rham cohomology of $C/k$, $\hdr^1(C/k)$,
which is the hypercohomology group of the de Rham complex $$ 0 \lra
\cO_C \lra \Omega_{C/k}^1 \lra 0$$ on $C$. We have a short exact
sequence
 \begin{equation}\label{2.1}
   0 \lra H^0(C, \Omega_{C/k}^1)
   \lra \hdr^1(C/k) \lra H^1(C, \cO_{C}) \lra 0,
\end{equation}
where, relying on the description of de Rham cohomology in terms of
forms of the second kind we have
\begin{itemize}\item $H^0(C,\Omega^1_{C/k})$, the space of
  differentials of the first kind, is identified with its image. It
  has dimension $g$, and we will denote it $\hdr^{1,0}(C/k)$.
 \item $H^1(C, \cO_{C})$ also has dimension $g$ and may be canonically
   identified with the tangent space at the origin of the Jacobian of
   $C$,  $J =\Pic^0(C)$.
 \item $\hdr^1(C/k)$ has a canonical non-degenerate alternating form
   given by the algebraic cup product
   pairing
   \begin{align*}
     \hdr^1(C/k) \times \hdr^1(C/k) &\lra k \\
     ([\mu_1],[\mu_2])\quad\quad\quad &\mapsto [\mu_1] \cup
     [\mu_2],\end{align*} which can be described by the formula
   \begin{equation}
    [\mu_1] \cup
     [\mu_2] = \sum_P \Res_P(\mu_2 \int \mu_1)\label{eq:Serre},
  \end{equation}where $\mu_1, \mu_2$ are differentials of the second kind, with
   classes $[\mu_1]$ and $[\mu_2]$, respectively, in $\hdr^1(C/k)$ and
   the sum is over all points in $C$. The residue does not depend on
   the choice of a particular local integral for $\mu_1$ because
   $\mu_2$ is of the second kind and has no residue at any point.
\end{itemize}

We will also need the theory of \emph{Coleman integration}. Details will be
discussed in Section~\ref{coleman} but for now it suffices to know
that for each meromorphic form $\omega$ on $C$ and to each $D \in \Div^0(C)$
, the theory allows us to define the integral $\int_D \omega\in k$. In particular, for two $k'$-rational points $P$ and $Q$,
where $k'$
is any finite extension of $k$, it allows us to define the integral
$\int_P^Q \omega\in k'$. This construction is such that as $Q$ (or
$P$) varies, this last integral is expressed as a power series in $Q$
which is locally convergent and whose differential is $\omega$. For
forms that have residues, the Coleman integral depends on the choice
of a branch of the $p$-adic logarithm function. We fix this choice for
the computation of the local height pairing to be the one determined
in~\eqref{logbranch}.

Let $T(k)$ denote the subgroup of differentials on $C$ of the third
kind. We have a \emph{residue divisor homomorphism}
\begin{equation*}
  \Res: T(k) \to \Div^0(C),\quad \quad \Res(\omega) = \sum_P \Res_P \omega
\end{equation*}
where the sum ranges over all closed points of $C$. That
the image is contained in $\Div^0(C)$ is just the residue theorem. By the residue divisor homomorphism, $T(k)$ fits into the
following exact sequence:
\begin{equation}
\xymatrix{0 \ar[r] & \Omega^1(C/k) \ar[r] & T(k) \ar[r]^{\Res}
  \ar[r] & \Div^0(C) \ar[r] & 0.}\label{someshort}
\end{equation}

 We are interested in a particular
subgroup of $T(k)$ whose elements are the logarithmic differentials,
i.e., those of the form $\frac{df}{f}$ for $f \in k(C)^*$. We denote
this subgroup as $T_l(k)$.  Since $T_l(k) \cap \hdr^{1,0}(C/k) =
\{0\}$ and $\Res(\frac{df}{f}) = (f)$, we deduce from the
sequence~\eqref{someshort} the short exact sequence
$$ 0 \lra
\hdr^{1,0}(C/k) \lra T(k)/T_l(k) \lra J(k) \lra 0.$$

This sequence has a natural identification with the $k$-rational points of an exact sequence of commutative algebraic groups over $k$:
$$0 \lra  H^0(C, \Omega_{C/k}^1) \lra E \lra J \lra 0,$$ where $E$ is the universal extension of $J$ by a vector group and $H^0(\Omega^1) \cong \Ga^g.$ Since the Lie algebra of $E$ is canonically isomorphic to $H_{dR}^1(C/k)$, the exact sequence (\ref{2.1}) is the resulting exact sequence of Lie algebras over $k$.

Now as $k$ is $p$-adic, we will make use of the fact that we have a
logarithmic homomorphism defined on an open subgroup of the points of
any commutative $p$-adic Lie group, $G$, to the points of its Lie
algebra $\Lie(G)$. When $G = E$ or $J$, the open subgroup on which the
logarithm converges has finite index, so the homomorphism can be
uniquely extended to the entire group. We denote this extension as
$\log_E$ or $\log_J$, respectively. Since the logarithm is functorial
and equal to the identity on $ H^0(C, \Omega_{C/k}^1)$, we obtain the
following:

\begin{proposition}\label{comm}There is a canonical
  homomorphism $$\Psi: T(k)/T_l(k) \lra \hdr^1(C/k)$$ which is the
  identity on differentials of the first kind and makes the following
  diagram commute:
  $$\xymatrix{
    0 \ar[r] & \hdr^{1,0}(C/k) \ar[r]\ar@{=}[d] & E(k) \ar[d]^{\Psi = \log_E}
    \ar[r] & J(k) \ar[r]\ar[d]^{\log_J}  & 0  \\
    0 \ar[r] & \hdr^{1,0}(C/k) \ar[r] &\hdr^1(C/k) \ar[r]
    &H^1(C,\cO_{C/k}) \ar[r] &0.}$$
\end{proposition}

Note that the map $\Psi$ takes a differential of the third kind on $C$
to a differential of the second kind modulo exact differentials. It
can be extended to a linear map from the $k$-vector space of all
differentials on $C/k$ to $\hdr^1(C/k)$ by writing an arbitrary
differential $\nu$ as a linear combination $\nu = \sum\alpha_i\mu_i +
\gamma$, where $\mu_i$ is of the third kind, $\alpha_i \in
\overline{k}$, and $\gamma$ is of the second kind on $C$. We then
define $\Psi(\nu) = \sum\alpha_i \Psi(\mu_i) + [\gamma]$.

The definition of the log map $\Psi$ is not very useful for
computations. An equivalent alternative definition has been given
in~\cite{besser:k2}. It is based on the notions of local and global
symbols or indices.
\begin{definition}\label{symboldef}
  For $\omega$ a meromorphic form and $\rho$ a form of
  the second kind, we define the \emph{global symbol} $\langle
  \omega,\rho\rangle$ as a sum of \emph{local symbols} $\langle
  \omega,\rho\rangle_A$. We have $$\langle \omega,\rho
    \rangle = \sum_{A} \langle \omega, \rho \rangle_A,$$ where \begin{equation}\label{symb2}\langle \omega, \rho \rangle_A = \Res_A\left(\omega\left(\int\rho +
        \int_Z^A\rho\right)\right).\end{equation} Note that the sum is taken over all points $A$ where either $\rho$ or $\omega$ has a singularity.  Each local
  symbol $\langle \omega, \rho \rangle_A$ is computed in the local
  coordinates of $A$, and $Z$ is a point that is fixed throughout the
  calculation of a single global
  symbol, which sets the constant of integration for the first indefinite integral of $\rho$. More precisely, we note that \eqref{symb2} is computed by rewriting $\omega$ and $\rho$ using the local coordinates at $A$, computing the indefinite integral of $\rho$ (as a formal power series, and in particular, not as a Coleman integral) and computing the Coleman integral of $\rho$ from $Z$ to $A$.
\end{definition}

The following result~\cite{besser:k2} reduces the computation of
$\Psi$ to the computation of global symbols.
 \begin{proposition}\label{2.6}
   Let $\omega$ be a meromorphic form
   and $\rho$ a form of the second kind. Then $\langle \omega, \rho
   \rangle = \Psi(\omega) \cup [\rho].$
\end{proposition}

\begin{proof}
See \cite[Prop 4.10]{besser:k2}.
\end{proof}

Now recall that we have at our disposal the complementary subspace
$W=W_v$. It allows us to isolate a canonical form $\omega_D$ with
residue divisor $D$ as follows.
\begin{definition}\label{omd}
  For any divisor ${D}$ of degree $0$ on $C$ we let
  $\omega_{D}$ be the unique form of the third kind
  satisfying $$\res(\omega_{D})={D}, \quad \lag(\omega_{D})\in W.$$
\end{definition}
It is easy to see from the properties of $\lag$ that this indeed uniquely defines the form
$\omega_{D}$.

\begin{definition}
  The local height pairing is defined by
  \begin{equation*}
    h_v({D_1},{D_2}):=t_v\left(\int_{D_2}\omega_{D_1}\right)
  \end{equation*}
  (recalling that the supports of ${D_1}$ and ${D_2}$ are
  disjoint), where $t_v$ is the trace map
  determined by the decomposition of $\ell_v$.
\end{definition}

\begin{remark}
In certain cases, there is a canonical complement $W$
to $\hdr^{1,0}(C/k)$ in $\hdr^1(C/k)$. Namely, when $C$ has good
ordinary reduction, we may take $W$ to be the unit root subspace for
the action of Frobenius.
\end{remark}

Some properties of the local height pairing are as follows:
\begin{proposition}\label{locprop}The local height pairing
  $h_v(D_1,D_2)$ is continuous and bi-additive. It is
  symmetric if and only if the subspace $W$ of $\hdr^1(C/k)$ is
  isotropic with respect to the cup product pairing. Finally, the
  formula~\eqref{eq:heighprinc} holds.
\end{proposition}

\begin{proof}
See \cite[Prop 5.2]{coleman-gross}.
\end{proof}

%% file: hyperelliptic.tex
\section{Hyperelliptic curves}
\label{sec:hyper}

In this section we review relevant facts about algebraic and
computational aspects of the theory of hyperelliptic curves.

Let us now suppose that $C$ is a hyperelliptic curve. As $k$ has characteristic not equal to 2, a hyperelliptic curve $C$ over $k$ of genus $g$
is an algebraic curve given by the equation \begin{equation}
  \label{eq:hyperell}
  y^2= f(x)\;,\; f(x) \in k[x],
\end{equation}
where $f$ has simple zeros. For our implementation in Sage, we shall further assume $f$ to be monic and $\deg f = 2g+1$.  

The curve $C$ is singular only at infinity (and non-singular when $g=1$).
To describe the neighborhood of infinity we normalize the curve there and obtain
the \emph{equation at infinity}
\begin{equation}
  \label{eq:infinity}
  t^2= s \frev(s) \text{ with } x=\frac{1}{s},\;\; y= \frac{t}{s^{g+1}},
\end{equation}
where $\frev(s) = s^{2g+1} f(1/s)$ is the reversed polynomial. The
point at infinity corresponds in these coordinates to both $s$ and $t$
being $0$. Furthermore, $t$ has a simple zero and $s$ has a double
zero at this point.

   As is well-known, the first de Rham cohomology of $C$ has a basis consisting of the forms of the second kind
\begin{equation*}
  \omega_i:= \frac{x^i dx}{2y} \text { for } i=0,\ldots, 2g-1.
\end{equation*}
We will denote this basis as \begin{equation}\label{eq:stanbasis}\mathcal{B} = \left\{ \frac{dx}{2y}, \frac{x dx}{2y}, \ldots, \frac{x^{2g-1}dx}{2y}\right\}.\end{equation}
If we make the change of coordinates \eqref{eq:infinity} we see that
these are transformed as follows:
\begin{equation*}
   \frac{x^i dx}{2y} \mapsto -s^{g-1-i} \frac{ds}{2t}.
\end{equation*}
Since $s$ has a double zero at the point at infinity one sees that
these forms are holomorphic for $i=0,\ldots,g-1$ and meromorphic
otherwise. 

We finally recall the \emph{hyperelliptic involution} $w$ defined by $w(x,y)=(x,-y)$.



%% file: coleman-integration.tex
\section{Coleman integrals}\label{coleman}
Here we review the relevant background on Coleman integrals and
describe new techniques to handle Coleman integrals of meromorphic
differentials with poles in non-Weierstrass residue discs.  This gives us the necessary tools to present our algorithm to compute local heights.
\subsection{Differentials of the second kind}
The foundational reference for this is \cite{coleman:torsion}. A more
expanded version of our presentation can be found in \cite{bbk}.

Let $\omega$ be a 1-form, with
$(\omega)_{\infty}$ denoting its polar support. For $P,Q \in C(k)$,
Coleman integration allows us to compute $\int_P^Q \omega \in k$.

\begin{theorem}\label{bigthm}Let $\mu, \nu$ be 1-forms on $C$ and
  $P,Q,R \in C(k)$. The (definite) Coleman integral has the
  following properties:
  \begin{enumerate}\item Linearity: $\int_P^Q (a\mu + b\nu) =
    a\int_P^Q \mu + b\int_P^Q \nu,$ for $P,Q \not\in (\mu)_{\infty}
    \cup (\nu)_{\infty}$.
  \item Additivity: $\int_P^R \mu = \int_P^Q \mu + \int_Q^R \mu$, for
    $P,Q,R \notin (\mu)_{\infty}$. \item Change of variables: If $C'$
    is another curve and $\phi: C \ra C'$ a rigid analytic map between
    wide opens then $\int_P^Q \phi^*\mu = \int_{\phi(P)}^{\phi(Q)}
    \mu$ if $\phi(P),\phi(Q) \not\in (\mu)_{\infty}$.\item Fundamental
    theorem of calculus: $\int_P^Q df = f(Q)-f(P)$ for $f$ a
    meromorphic function $f$ on a wide open
    subset. \end{enumerate}\end{theorem}
\begin{proof} See \cite[Thm 2.3, Prop 2.4, Thm 2.7]{coleman:torsion}
  for details.\end{proof} Thus writing $$\omega = df + c_0\omega_0 +
\cdots + c_{2g-1}\omega_{2g-1},$$ where $f$ is a function and $c_i \in
k$, we see that the problem of computing $\int_P^Q \omega$ reduces
to computing the Coleman integral of a basis differential, $\int_P^Q
\omega_i$.

For hyperelliptic curves, Coleman integration can be performed
numerically. For simplicity we describe things only over $\Q_p$
although things work in complete generality (see~\cite{bbk} and also
\cite{Bes10a} where a slightly different but related version, which
was implemented in~\cite{Gut06}, is described).

Let us recall the following algorithms (Algorithms 8 and 11,
respectively) from \cite{bbk}. Throughout, let $\phi$ denote a $p$-power lift
of Frobenius.

\begin{algorithm}[Tiny Coleman integrals]\label{tiny} \;\quad\;\\
Input: \begin{itemize}\item Points $P, Q \in C(\Q_p)$ in the same residue disc (neither equal to the point at infinity).
\item A basis differential $\omega_i$.\end{itemize}
Output: The integral $\int_P^Q \omega_i$.\\
The algorithm:
\begin{enumerate}
\item Construct a linear interpolation from $P$ to $Q$. For instance, in a
non-Weierstrass residue disc, we may take
\begin{align*} x(t) &= (1-t)x(P) + tx(Q)\\
y(t) &= \sqrt{f(x(t))},\end{align*} where $y(t)$ is expanded as a formal power series in $t$.

\item Formally integrate the power series in $t$:$$\int_P^Q \omega_i = \int_P^Q x^i\frac{dx}{2y} = \int_0^1 \frac{x(t)^i}{2y(t)}\frac{dx(t)}{dt}dt.$$\end{enumerate}\end{algorithm}

\begin{algorithm}[Coleman integration in non-Weierstrass discs]\label{algo:nonteich}\;\quad\;\\
Input:\begin{itemize}\item The basis differentials $\{\omega_i\}_{i=0}^{2g-1}$.
\item Points $P,Q \in C(\Q_p)$ in non-Weierstrass residue discs.\end{itemize}
Output: The integrals $\{ \int_P^Q \omega_i\}_{i=0}^{2g-1}$.\\
The algorithm:
\begin{enumerate}
\item
Calculate the action of Frobenius on each basis element :
$$\phi^*(\omega_i) = df_i + \sum_{j=0}^{2g-1} M_{ij}\omega_j.$$

\item
By change of variables
(see Remark~\ref{rem:nonteich}), we obtain
\begin{equation}\label{linear}\sum_{j=0}^{2g-1}
  (M-I)_{ij}\int_P^Q \omega_j = f_i(P)-f_i(Q) -
\int_P^{\phi(P)}\omega_i -\int_{\phi(Q)}^Q\omega_i.\end{equation}
As the eigenvalues of the matrix $M$ are algebraic integers of
$\C_p$-norm $p^{1/2} \neq 1$ (see \cite[\S 2]{kedlaya:mw}), the matrix $M-I$ is invertible,
and we may solve \eqref{linear} to obtain the integrals $\int_P^Q \omega_i$.

\end{enumerate}
\end{algorithm}

\begin{remark} \label{rem:nonteich}
We obtain \eqref{linear} as follows.
By change of variables,\begin{align*}\int_{\phi(P)}^{\phi(Q)}\omega_i &= \int_{P}^{Q}\phi^*(\omega_i)\\
&= \int_{P}^{Q}(df_i + \sum_{j=0}^{2g-1}M_{ij}\omega_j)\\
&= f_i(Q)-f_i(P) + \sum_{j=0}^{2g-1}M_{ij}\int_P^Q \omega_j.\end{align*}
Adding $\int_{P}^{\phi(P)}\omega_i + \int_{\phi(Q)}^Q\omega_i$ to both sides of this equation yields $$\int_P^Q\omega_i = \int_P^{\phi(P)}\omega_i + \int_{\phi(Q)}^Q\omega_i + f_i(Q) - f_i(P) + \sum_{j=0}^{2g-1}M_{ij}\int_P^Q \omega_j,$$ which is equivalent to \eqref{linear}.
\end{remark}

\subsection{Meromorphic differentials}\label{mer}
The above approach does not work for a meromorphic
differential with poles in non-Weierstrass residue discs. We provide a new approach for dealing with this case (again we describe things only over $\Q_p$ but things work in general). Let $\omega$ be such a differential. As before, if $R,S$ are points in the same non-Weierstrass residue disc (different from those discs containing the poles of $\omega$), then $\int_S^R \omega$ is just a tiny integral, which can be computed as in Algorithm~\ref{tiny}. Let us now suppose that $R,S$ are in different non-Weierstrass residue discs.
\begin{algorithm}[Coleman integration: differential with poles in non-Weierstrass discs]\label{algo:mer}\;\quad\;\\
Input:\begin{itemize}\item The differential $\omega$ with residue divisor $(P)-(Q)$, with non-Weierstrass points $P,Q \in C(\Q_p)$.
 \item Points $R,S \in C(\Q_p)$ in distinct non-Weierstrass residue discs.\end{itemize}
Output: The integral $\int_S^R \omega$.\\
The algorithm:
\begin{enumerate}
\item Let $\alpha = \phi^*(\omega) - p\omega$. Using the methods of Section~\ref{subs:psi}, compute $\Psi(\omega)$. Using this, compute $\Psi(\alpha)$ as $\phi^*(\Psi(\omega)) - p\Psi(\omega)$.
\item Let $\beta$ be a form with residue divisor $(R)-(S)$. Compute $\Psi(\beta)$.
\item Compute the cup product $\Psi(\alpha) \cup \Psi(\beta)$ (see Section~\ref{subs:cups} for more details).
\item Evaluate the tiny integrals $\int_{\phi(S)}^S \omega$ and $\int_{R}^{\phi(R)} \omega$.
\item\label{alg3:res} Let $\mathcal{S}$ be the set of Weierstrass points and poles of $\alpha$. For each point $A \in \mathcal{S}$, expand $\alpha \int\beta$ in terms of the local coordinate at $A$ and find the residue at $A$. Compute the sum $$\sum_{A\in\mathcal{S}} \Res_A(\alpha\int\beta).$$
\item By Coleman reciprocity (see Remark~\ref{rem:recip}), recover the desired integral: \begin{equation}\label{alpha}\int_S^R \omega =
  \frac{1}{1-p}\left(\Psi(\alpha) \cup \Psi(\beta) + \sum_{A\in\mathcal{S}} \Res_A
    \left(\alpha\int\beta\right) - \int_{\phi(S)}^S \omega -
    \int_{R}^{\phi(R)}\omega\right).\end{equation}
\end{enumerate}
\end{algorithm}

\begin{remark}\label{rem:recip}
We obtain \eqref{alpha} as follows: \begin{align}\label{int}\int_{S}^{R} \alpha &= \int_S^R \phi^*\omega - p\int_S^R \omega \\
  &= \int_{\phi(S)}^{\phi(R)} \omega -p\int_S^R \omega \nonumber \\
  &= (1-p)\int_{S}^{R} \omega + \left(\int_{\phi(S)}^S \omega +
    \int_R^{\phi(R)}\omega \right)\nonumber.\end{align} As $\beta$ is a
form with residue divisor $(R)-(S)$, by Coleman reciprocity
\cite[Theorem 5.2]{coleman:reciprocity}, we compute the LHS of (\ref{int}) to be
$$\sum_{A \in \mathcal{T}} \Res_A\left(\beta\int\alpha\right) =  \Psi(\alpha) \cup \Psi(\beta)
- \sum_{A \in \mathcal{S}}
\Res_A\left(\alpha\int\beta\right),$$
where $\mathcal{T}$ is the set of Weierstrass points and poles of $\beta.$ We then have
\begin{equation*}\int_S^R \omega =
  \frac{1}{1-p}\left(\Psi(\alpha) \cup \Psi(\beta) + \sum_{A\in \mathcal{S}} \Res_A
    \left(\alpha\int\beta\right) - \int_{\phi(S)}^S \omega -
    \int_{R}^{\phi(R)}\omega\right).\end{equation*}
\end{remark}

\begin{remark}
In practice, the computation in Step~\ref{alg3:res} of Algorithm~\ref{algo:mer} is the slowest part of the algorithm, as it involves high-precision local calculations over all poles of $\alpha$ and all Weierstrass points of the curve. However, since $\sum_{T\in U}\Res_T(\alpha)= 0$ in each residue disc $U$, for the Weierstrass discs, we do not need a constant of integration. For the poles of $\alpha$, we may choose one constant of integration within each residue disc. More precisely, if $P$ and $Q$ are in separate residue discs, we compute
\begin{align*}\sum_{A \in U_P} \Res_A\left(\alpha\int\beta\right) &= \tr_{k(x(P_1))/k}\left(\Res_{P_1}\left(\alpha\int_P^{P_1}\beta\right)\right) = \tr_{k(x(P_1))/k}\left(\int_P^{P_1}\beta\right) ,\\
  \sum_{A \in U_Q}\Res_A\left(\alpha\int\beta\right) &=
  \tr_{k(x(Q_1))/k}\left(\Res_{Q_1}\left(\alpha\int_Q^{Q_1}\beta\right)\right)
  = \tr_{k(x(Q_1))/k}\left(-\int_Q^{Q_1}\beta\right) ,\end{align*}
where $P_1$ (resp, $Q_1$) is a root of $\alpha$ in the residue disc of
$P$ (resp $Q$).\end{remark} 

%% file: algorithm.tex
\section{The local height pairing at primes above $p$}
\label{sec:localatp}

In this section we will explain the algorithm that computes the local height
pairing at a prime above $p$ for degree zero divisors on the hyperelliptic curve $C$. Recall that we have as additional data the complementary subspace
$W$ and the character $\chi$ from which we deduce a branch of the
logarithm to be used in Coleman integration and the trace map $t_v$
(we keep the subscript $v$ at some places for clarity, even though it
now serves no purpose).

Let $D_1$ and $D_2$ be two divisors of degree $0$ on $C$.
Our main algorithm computes the local height pairing $h_v(D_1,D_2)$. It may be
described in two steps
\begin{itemize}
\item Compute the height pairing in the case where $D_1$ and $D_2$ are
  anti-symmetric  with respect to
  the hyperelliptic involution (Algorithm~\ref{antisymmetricalgorithm})
\item Compute the height pairing in the general case using the first
  case (Algorithm~\ref{generalalgorithm}).
\end{itemize}

Before discussing either algorithm, we begin with some general notes
about the representation of divisors on hyperelliptic curves
(see~\cite[App \S 5-6]{Kob98}).

Recall that a divisor of degree
$0$ on $C$ may be written in the form $$D=\sum m_i P_i - \left(\sum m_i\right) (\infty).$$

\begin{definition}
  A divisor $D$ as above is called \emph{semi-reduced} if the following conditions
  are satisfied:
  \begin{itemize}
  \item $m_i\ge 0$
  \item If $P_i$ is in the support of $D$, then $-P_i$ is not, unless
    $P_i=-P_i$ in which case $m_i=1$.
  \end{itemize}
  A semi-reduced divisor is called \emph{reduced} if in addition
  \begin{itemize}
  \item $\sum m_i \le g$.
  \end{itemize}
\end{definition}

One may represent a semi-reduced divisor $D$
by a pair of
polynomials $a(x),b(x)$ with $\deg(b)< \deg(a)$ such that
\begin{itemize}
\item The projection of $\sum m_i P_i$ on $\PP^1$ is the zero divisor of $a(x)$.
\item $b(x)$ is an interpolation polynomial with the property that for
  $P_i= (x_i,y_i)$ we have $b(x_i)=y_i$.
\end{itemize}
The condition that $D$ is
reduced is equivalent to having  $\deg(a)\le g$.
\begin{remark}
  One can associate $b$ uniquely to the divisor by insisting that
  $a(x) | (b(x)^2-f(x)) $. This would be less important for us and there
  are cases we may not achieve this.
\end{remark}
\begin{definition}
  Let us denote by $(a,b)$ the semi-reduced divisor determined by the
  pair of polynomials $a$ and $b$ and call $(a,b)$ the standard
  representation of the divisor.
\end{definition}

It is known that any degree zero divisor $D$ on $C$ is equivalent to a
unique reduced divisor. Furthermore, the reduction is
effective. More precisely, passing from an arbitrary divisor to a
semi-reduced one is just a question of adding or subtracting divisors
of functions pulled back from $\PP^1$ while passing from a semi-reduced
divisor to a reduced divisor has an effective algorithm described
in~\cite[App, Alg~2 and Thm~7.2 ]{Kob98}. Since we know
our height pairing satisfies~\eqref{eq:heighprinc} by
Proposition~\ref{locprop}, which easily
allows to pass from a divisor to an equivalent divisor in the pairing
there is no harm in  assuming that our divisors are reduced.

Unfortunately for us, reduced divisors are not sufficient. The reason
is that since they always have a component at infinity, two such
divisors cannot have disjoint support unless one of them is
trivial. For this reason we will work with the difference of two
reduced divisors.
\begin{definition}
  The divisor denoted $(a,b)-(c,d)$, where $a$ and $c$
  are polynomials of the same degree $\le g$, stands for the
  difference of the reduced divisors defined by $(a,b)$ and $(c,d)$.
\end{definition}
We always assume that the two divisors defined by $(a,b)$ and $(c,d)$
have no common components. If there are common components they can be
cancelled out.

We will mostly work with antisymmetric divisors. Given any zero
divisor $D$, the divisor $D-w^\ast D$ is antisymmetric. Conversely, any
antisymmetric divisor is obtained in this way. It follows easily that
any antisymmetric divisor is equivalent to $D-w^\ast D$ for a reduced
divisor. There may be several representations in this form for a given
divisor, however, there is just one containing no points $P_i$ with
$P_i=-P_i$. In the representation $(a,b)$ for $D$ this is equivalent to having
$a$ prime to $f$.
\begin{definition}\label{antirep}
  The standard representation of an anti-symmetric divisor is
  $$[a,b]:=D-w^\ast D$$
  with $D$ a reduced divisor given in the form $(a,b)$. It is a
  reduced standard representation if $a$ is prime
  to $f$.
\end{definition}

Note that with anti-symmetric divisors the standard representation
suffices because they do not have components at infinity. Thus,
$[a_1,b_1]$ and $[a_2,b_2]$ have disjoint support if and only if
$a_1$ and $a_2$ are relatively prime.

We now describe Algorithm~\ref{generalalgorithm}. For any divisor $D$ we have a
decomposition, with rational coefficients
\begin{equation}
  \label{eq:decomp}
  D= \frac{1}{2} D^+ +  \frac{1}{2} D^-, \quad \quad D^+ = D+
  w^\ast(D)\;,\;  D^- = D- w^\ast(D).
\end{equation}
\begin{lemma}\label{minuspart}
  Suppose $D$ is given by the representation $(a,b)-(c,d)$. Then, the
  divisor $D^-$ is just $[a\cdot c, e]$ in terms of
  Definition~\ref{antirep}, where $e$ is obtained by solving a Chinese
  Remainder Theorem problem to be congruent to $b$ modulo $a$ and to
  $-d$ modulo $c$.
\end{lemma}
\begin{proof}
  If $a$ is prime to $c$ this is clear. In general, suppose
  $(x-\alpha)$ has multiplicity $m$ in $a$ and $n$ in $c$. We may
  assume $m\ge n$. Since we are assuming $(a,b)$ and $(c,d)$ have no
  common components, it follows that $D-w^\ast(D)$ is going to have the two
  summands $(m+n)[(\alpha,b(\alpha))-(\alpha,-b(\alpha)]$, and $m+n$
  is indeed the multiplicity of $ac$ in $\alpha$ (it seems though that
  we can only solve the Chinese remainder problem modulo the least
  common multiple of $a$ and $c$).
\end{proof}
On the other
hand, $D^+$
is nothing but the divisor of the rational function $\frac{a(x)}{c(x)}$
considered as a rational function on $C$. It follows from the fact
that~\eqref{eq:heighprinc} is satisfied by Proposition~\ref{locprop} that for
any $E\in \Div^0(C)$ we have
\begin{equation}
  h_v(D^+,E)= h_v(E,D^+)=\chi\left(\frac{a}{c}(E)\right),\label{eq:withsymm}
\end{equation}

 where $\frac{a}{c}(E)$ means as usual the product of
the values of $\frac{a}{c}$ on the $x$-coordinates of the points making up $E$
with the appropriate multiplicities. An easy consequence of this
formula is that
\begin{equation}\label{symmantisymm}
  h_v(D^+,E)= h_v(E,D^+)=0 \quad \text{if } E \text{ is antisymmetric}.\end{equation}

Consider now two divisors $D_1$ and $D_2$ in $\Div^0(C) $. Decomposing
into plus and minus parts it follows from~\eqref{symmantisymm} that
\begin{equation}\label{generalexpression}
  h_v(D_1,D_2) =  \frac{1}{4} h_v(D_1^+,D_2^+)+ \frac{1}{4}  h_v(D_1^-,D_2^-).\end{equation}
The first term can be computed
using~\eqref{eq:withsymm}, while the second term is a height pairing
between anti-symmetric divisors. This immediately gives the following
algorithm.
\begin{algorithm}[$p$-adic height pairing for general
  divisors]\label{generalalgorithm}\;\quad\;\\
  Input:
  \begin{itemize}
  \item The subspace $W$, branch of logarithm and trace map $t$.
  \item Divisors $D_1$ and $D_2$ with disjoint support given as
    \begin{align*}D_1 &= (a_1,b_1)-(c_1,d_1) \\
    D_2 &=  (a_2,b_2)-(c_2,d_2).\end{align*}
  \end{itemize}
  Output: The local height pairing $h_v(D_1,D_2)$.\\
  The algorithm:
  \begin{enumerate}
  \item Compute expressions for the divisors $D_1^-$ and $D_2^- $ using Lemma~\ref{minuspart}.

  \item Compute using~\eqref{eq:withsymm},
    \begin{equation*}
       h_v(D_1^+,D_2^+) = \chi((a_1/c_1)(D_2^+)).
    \end{equation*}
  \item Compute, using Algorithm~\ref{antisymmetricalgorithm}, the local
    height pairing for anti-symmetric divisors $ h_v(D_1^-,D_2^-) $.
  \item Substitute in~\eqref{generalexpression} to obtain $h_v(D_1,D_2)$.
  \end{enumerate}
\end{algorithm}

We next turn to Algorithm~\ref{antisymmetricalgorithm} for the case of
anti-symmetric divisors $D_1$ and $D_2$. First of all, we have to introduce yet
another decomposition. The algorithm behaves differently with respect
to parts that reduce to the Weierstrass points (Weierstrass divisors) and those which do
not. We can decompose a divisor $D$ into the sum of its Weierstrass
part $D\Rw$ and its non-Weierstrass part $D\Rnw$. Then, in a similar
way to~\eqref{generalexpression} we have the decomposition
\begin{equation}\label{weierdecomp}
  h_v(D_1,D_2) =  h_v(D_1,D_2\Rw)+
  h_v(D_1\Rw,D_2\Rnw) +   h_v(D_1\Rnw,D_2\Rnw).
\end{equation}

We now give the algorithm and later we discuss
each step with some more details.
\begin{algorithm}[$p$-adic height pairing for anti-symmetric
  divisors]\label{antisymmetricalgorithm}\;\quad\;\\
  Input:
  \begin{itemize}
  \item The subspace $W$, branch of logarithm and trace map $t$.
  \item Anti-symmetric divisors $D_1$ and $D_2$ given in standard
    representation $(a_1,b_1)$, $(a_2,b_2)$.
  \end{itemize}
  Output: The local height pairing $h_v(D_1,D_2)$.\\
  The algorithm:
  \begin{enumerate}
  \item Compute the cup product matrix for a basis of $\hdr^1(C/k)$, as in Section~\ref{subs:cups}.
  \item Compute  $D_1\Rw$ and $D_1\Rnw$.
  \item Write down forms $\nu_1\Rw$ and $\nu_1\Rnw$ with
    residue divisors $D_1\Rw$ and $D_1\Rnw$, respectively.
  \item Compute the form $\omega_{D_1\Rw}$ and a holomorphic form
    $\eta$ such that $$\omega_{D_1\Rnw}=\nu_1\Rnw -\eta.$$
  \item Compute the tiny  Coleman integral $ h_v(D_1,D_2\Rw) =t( \int_{D_2\Rw}
    (\omega_{D_1\Rw}+ \nu_1\Rnw -\eta))$.
  \item Compute the Coleman integral $ h_v(D_1\Rw,D_2\Rnw)= t
    (\int_{D_2\Rnw} \omega_{D_1\Rw}) $.
  \item Compute the Coleman integral $ h_v(D_1\Rnw,D_2\Rnw) = t
    ( \int_{D_2\Rnw} \nu_1\Rnw -  \int_{D_2\Rnw} \eta)$.
  \item Compute $h_v(D_1,D_2)$ using~\eqref{weierdecomp}.
  \end{enumerate}
\end{algorithm}

We now add some further details on each step in this algorithm.


\subsection{Computing cup products}\label{subs:cups}

We first compute the cup product between any two elements of the
standard basis~\eqref{eq:stanbasis} for
$\hdr^1(X/k)$. This is easily done using the formula~\eqref{eq:Serre}.

We can be a bit more precise as follows.
\begin{definition}The \emph{cup product matrix}
associated to $C$  with respect to $\mathcal{B}$ is the $2g \times 2g$ matrix with entry $a_{i,j}$ given by the cup product of differentials
$[\omega_{i-1}]\cup[\omega_{j-1}]$, normalized so that $[\omega_{i-1}]\cup [\omega_{j-1}] = \Res(\omega_{j-1}\int\omega_{i-1})$ .
\end{definition}

By computing in the local coordinates at infinity, we may record the following
\begin{lemma}The cup product matrix for $C$ with respect to
  $\mathcal{B}$ satisfies the following properties:
  \begin{enumerate}\item
    Anti-diagonal elements are given by the sequence
$$\left\{\frac{1}{2g-1},\frac{1}{2g-3},\ldots,\frac{1}{3}, 1, -1,
  -\frac{1}{3},\ldots,-\frac{1}{2g-1}\right\}.$$
\item Entries above the anti-diagonal are 0. \item Diagonal elements
  are 0.
\end{enumerate}\end{lemma}
\begin{example}\label{cpm-gen1}The cup product matrix with respect to
  $\mathcal{B}$ for an elliptic curve  is $\left( \begin{array}{cc} 0
      & 1  \\
      -1 & 0  \end{array} \right).$ Note that the subspace spanned by
  $\frac{xdx}{2y}$ is isotropic. In particular, for genus 1, we may
  take $\hdr^{1,0}(C/k)$ spanned by $\omega_0$ and $W$ spanned by
  $\omega_1$, and we need not require $p$ to be a prime of ordinary
  reduction for the pairing to be symmetric.\end{example}

\subsection{The map $\Psi$}\label{subs:psi}
We compute $\Psi$ of a differential $\omega$ by writing $$\Psi(\omega) = c_0\omega_0 + \cdots + c_{2g-1}\omega_{2g-1},$$ and solving for the coefficients $c_i$. This is done by considering a linear system involving global symbols and cup products: \begin{align*}\langle
\omega,\omega_j\rangle &= \Psi(\omega) \cup [\omega_j] = \sum_{i=0}^{2g-1} c_i([\omega_i]\cup
[\omega_j]).\end{align*}

Recall that as in Definition~\ref{symboldef}, we calculate the global symbol as a
sum of local symbols, each of which involves a Coleman integral and a
calculation in local coordinates. This computation is actually much
simpler:
\begin{proposition}
  Suppose $\omega$ is a form of the third kind with residue divisor
  $D$ which does not contain
  $\infty$. Then we have
  \begin{equation*}
    \langle \omega, \omega_i \rangle = \int_D \omega_i + \Res_\infty\left(\omega\int \omega_i\right),  \end{equation*}
  where the residue at $\infty$ is computed by taking any antiderivative of   $\omega_i$.\end{proposition}

\begin{proof}
The sum of local symbols is over all points where either
$\omega$ or $\omega_i$ has a singularity. These are the points in the support of $D$ and possibly the point $\infty$. Since $\omega$ has
a simple pole at each point $P$ in the support of $D$, the local
symbol is simply the multiplicity of $D$ at $P$ times $\int_Z^P
\omega_i$ (where $Z$ is a fixed point throughout the global symbol calculation). Summing over all points gives $\int_D \omega_i$. On the
other hand, since we are assuming that $\omega$ is holomorphic at
$\infty$, the choice of the constant of integration for $\omega_i$ at
$\infty$ does not matter.
\end{proof}

Now letting $N$ denote the cup product matrix, we have $$\Psi(\omega) =  N\inv\left( \begin{array}{c}
 \langle \omega, \omega_0 \rangle \\
 \vdots \\
\langle \omega, \omega_{2g-1} \rangle  \end{array} \right).$$

\subsection{Decomposing a divisor $D$ into $D\Rw$ and $D\Rnw$}
\label{sec:weierdecomp}

This is very easy to do. When the divisor is given in standard form
$[a,b]$ one just reduces $a$ modulo the prime $\pi$, picks up the part
that reduces to Weierstrass points by taking the greatest common
divisor with the reduction of $f$, and then applies a Hensel lift to get
the factor $a\Rw $ of $a$ corresponding to the points reducing to
Weierstrass. Then we have $a\Rnw = a/ a\Rw$ and the divisor
decomposition is deduced from that.

\subsection{A form with the required residue divisor}
\label{sec:formresidue}

This is an easy task with the following:
\begin{proposition}\label{getform}
  Let the anti-symmetric divisor $D$ be given in standard
  representation $[a,b]$. The differential form
  \begin{equation*}
    \omega = \frac{a'(x) b(x)dx}{a(x) y}
  \end{equation*}
  has simple poles and its residue divisor is $D$.
\end{proposition}
\begin{proof}
Suppose that $(a,b)=\sum m_i P_i - \left(\sum m_i \right)(\infty)$. We can write $\omega=
(b/y) \dlog(a)$. The form $\dlog(a)$ has
simple poles at $\pm P_i$ with residue $m_i$ while $b/y$ has value $1$
at $P_i$ and value $-1$ at $-P_i$. On the other hand we can also write
$\omega= (a' b/a) (\frac{dx}{2y})$ and since $\frac{dx}{2y}$ is holomorphic, it follows
that $\omega$ has no poles where $a$ does not vanish. Finally, it is
easy to see that $\omega$ does not have a pole at infinity.
\end{proof}

\subsection{Finding $\omega_D$ for a Weierstrass divisor $D$}

Suppose we have already written down a form whose residue divisor is
$D$. Since the singularities of the form are contained in the
Weierstrass residue discs, it is amenable to the reduction algorithm
done in Kedlaya's algorithm~\cite{kedlaya:mw}. This means that we may compute a
representation of $\omega$ as a linear combination of the basis
$\mathcal{B}$ plus an exact differential $dg$. Since it follows
that $\Psi(\omega) $ is just the above combination of basis elements, we need only subtract the appropriate combination of holomorphic
basis elements to make it reside in $W$.

\subsection{Finding $\omega_D$ for a non-Weierstrass divisor $D$}
\label{sec:omy}
We start with a form $\omega$ whose residue divisor is $D$. We compute $\Psi(\omega)$ as in Section~\ref{subs:psi}. All that is left to do is to let $\eta$ be the projection of $\Psi(\omega)$ on $\hdr^{1,0}$ along $W$.

\subsection{Integration when $D_2$ reduces to Weierstrass points}
\label{sec:ztoweir}

Suppose that the points of $D_2$ reduce to the Weierstrass
points. Since $D_2$ is anti-symmetric, this
means that in computing the integral we need to take the sum of
differences over pairs of points $\pm P $ which reduce to the same
point. This is a sum of tiny integrals.

\subsection{Computation when $D_1$ is Weierstrass and $D_2$ is not}
\label{sec:ztoweir}

In this case, we are given, by the previous reduction, the form
$\omega_{D_1}$ as a combination $\sum \alpha_i \omega_i +dg$. Thus
$\int_{D_2} \omega_{D_1}=\sum \alpha_i \int_{D_2} \omega_i +g(D_2)$.
Since the points of $D_2$ are in the domain where the integrals of the
$\omega_i$ may be computed, this is a standard computation.

\subsection{Computation when both divisors are non-Weierstrass}
\label{sec:away}

In this case, $\omega_{D_1}$ is given in the form $\nu_1 - \eta$,
where $\nu_1$ is a form with residue divisor $D_1$
and $\eta$ is holomorphic. The integral of $\eta$ poses no
problems while that of $\nu_1$ is discussed in Subsection~\ref{mer}. 

%% file: implementation.tex
\section{Implementation notes}\label{notes}

In this section, we discuss the choices made in our Sage~\cite{sage}
implementation and give error estimates on the precision of our
results.  We work over $\Q_p$ with a precision of $n$ digits; note that if one
desires an answer with $n$ digits of precision, one has to start with
a larger working precision, as seen below. 

We only discuss the
computation for anti-symmetric divisors, as the extension to general
divisors is trivial, as discussed in Section~\ref{sec:localatp}. Furthermore, our
implementation assumes that the divisors are of the form $(P)-(-P)$
for a $\Q_p$-rational point $P$, as it is then quite easy to consider cases
when the divisor is a sum of such expressions. Finally, all
computations are done with respect to a particular choice of the complementary subspace $W$, which we
describe below.

\subsection{Auxiliary choices}
Our algorithm relies on the splitting $$\hdr^1(C/k) \simeq H_{dR}^{1,0}(C/k)\oplus W,$$ which allows us to write $$\log(\nu_1) = \eta + \log(\omega_{D_1}),$$ where $\eta$ is holomorphic, and $\log(\omega_{D_1}) \in W$.

As noted in Example \ref{cpm-gen1}, when $g = 1$, we may simply take $\hdr^{1,0}(C/k)$ spanned by $\omega_0$ and $W$ spanned by $\omega_1$. Note that for genus $g >1$ it does not suffice to take $W$ spanned by $\omega_g,\omega_{g+1},\ldots,\omega_{2g-1}$, as the resulting subspace is not necessarily isotropic.  While the local height is independent of basis, it is not independent of the choice of $W$. For the local height to be symmetric, it is necessary that $W$ be an isotropic subspace. For $g> 1$, we thus further require that $p$ be a prime of \emph{ordinary} reduction, so that $W$ can be chosen to be the unit root subspace.
Generalizing \cite{mazur-stein-tate}, we compute a basis for $W$ as follows:

\begin{proposition}Let $n$ be the working precision in the underlying base ring $\Z_p$, so that all computations are done modulo $p^n$. Let $\Frob$ denote the matrix of a $p$-power lift of Frobenius, as acting on the standard basis $\mathcal{B}$ of $\hdr^1(C/k)$. Then $\left\{\Frob^n \omega_g, \Frob^n \omega_{g+1},\ldots, \Frob^n \omega_{2g-1}\right\}$ is a basis for $W$.\end{proposition}
\begin{proof}With the integral structure provided by crystalline cohomology, it is
well known that $\Frob$ maps the holomorphic forms to $p$ times the integral
structure. Thus, with $W$ the unit root part decomposing a vector $v$ into
$\omega+\eta$ with $\omega\in W$ and $\eta$ holomorphic it is easy to see that
$\Frob^n v$ is in $W+ p^n$ times the integral structure. In other words, up
to the prescribed precision, $\Frob^n v$ lies in $W$. On the other hand, $\Frob$ is
invertible so starting with $g$ independent vectors modulo the holomorphic differentials one gets $g$ independent vectors in $W$.\end{proof}

\subsection{Precision}\label{9.1}
Broadly speaking, the $p$-adic precision of a local height depends on
two types of calculation: \begin{enumerate}\item Coleman integrals of
  basis differentials (or otherwise ``nice'' differentials -- e.g.,
  holomorphic in the discs corresponding to the limits of integration)
  and \item expansion of local coordinates at a point.\end{enumerate}
Each key step of the algorithm in Section~\ref{algo:mer} can be categorized
as depending on one or both of these:
\begin{itemize}\item $\Psi(\omega)$ needs the cup product matrix
  (local coordinates) and local symbols (Coleman integrals of basis
  differentials)
\item $\int \eta$ is a sum of Coleman integrals of basis differentials
\item $\int \omega$ is defined in terms of \begin{itemize}\item tiny
    integrals (Coleman integrals of a ``nice'' differential), \item
    sums of residues of Laurent series (local coordinates), and \item
    $\Psi(\alpha)$, $\Psi(\beta)$ (as
    above).\end{itemize}\end{itemize}

\subsubsection{Precision of Coleman integrals}
For more details on the rigorous computation of Coleman integrals, see
\cite{bbk}. We recall the following results (\cite[Props 18-19]{bbk}):

\begin{proposition} \label{prop:tiny precision}
Let $\int_P^Q \omega$ be a tiny integral in a non-Weierstrass residue disc, with $P,Q$ defined over
an unramified extension of $\Q_p$ and accurate to $n$ digits of precision.
Let $(x(t),y(t))$ be the local interpolation between $P$ and $Q$ defined by
\begin{align*}x(t) &= x(P)(1-t) + x(Q)t = x(P) + t(x(Q)-x(P))\\
y(t) &= \sqrt{f(x(t))}.\end{align*}
Let $\omega = g(x,y)dx$ be a differential of the second kind such that
$h(t) = g(x(t), y(t))$ belongs to $\cO[[t]]$.
If we truncate $h(t)$ modulo $t^m$, then
the computed value of the integral $\int_P^Q \omega$ will be correct to
$\min\{n, m+1- \lfloor \log_p (m+1) \rfloor\}$ digits of (absolute) precision.
\end{proposition}

\begin{proposition}\label{m}Let $\int_P^Q \omega$ be a Coleman integral, with $\omega$ a differential of the second kind and with $P,Q$ in non-Weierstrass residue discs, defined over an unramified extension of $\Q_p$, and accurate to $n$ digits of precision.
Let $\Frob$ be the matrix 
of the action of Frobenius on the basis differentials.
Set $B = \Frob^t - I$, and let $m = v_p(\det(B))$. Suppose the relevant tiny integrals have series expansions truncated modulo $t^{n-1}$. 
Then the computed value of the integral $\int_P^Q \omega$ will
be accurate to $n-\max\{m , \lfloor \log_p n \rfloor\}$ digits of precision.\end{proposition}

We will now carefully review the precision of each of the objects we
computed, as an expansion of the overview in Section \ref{9.1}. Let
$\omega$ be a differential with residue divisor $$D_1 = (P) - (-P)$$ and $\beta$ a
differential with residue divisor $$D_2 = (Q) - (-Q).$$ The precision of
$\Psi(\omega)$ (and $\Psi(\beta)$) just depends on the Coleman
integral involved, as the residue can just be read off of the
differential.

After computing $\Psi(\omega)$ with respect to the standard basis of
$\hdr^1(C)$, we fix a splitting of $\hdr^1(C/k)$ into
$\hdr^{1,0}(C/k) \oplus W$, which gives $\eta$ and $\omega_{D_1} =
\omega-\eta$. Since the height pairing is given by $\int\omega_{D_1}$,
we need to compute the integrals $\int\omega$ and $\int\eta$.

The integral $\int\eta$ is just a linear combination of the integrals
of holomorphic basis differentials. On the other hand, the integral of
$\omega$ requires the computation of $\Psi(\alpha)$, $\Psi(\beta)$,
$\sum\Res(\alpha\int\beta)$, and the tiny integral $\int_Q^{\phi(-Q)}\omega$. As before, the tiny integral
is computed with precision as above.

Since $\alpha = \phi^*\omega-p\omega$, we may write $\Psi(\alpha)$ in
terms of things we have already computed, namely $\Psi(\alpha) =
\Frob(\Psi(\omega))-p\Psi(\omega)$. So need not do more work here.
However, the precision of $\sum\Res(\alpha\int\beta)$ merits further
discussion, as we must consider its representation in local
coordinates.

\subsubsection{Precision of local coordinates}
Computing with local coordinates is crucial to the algorithm. More
precisely, for any point $P$, we must construct power series $x(t),
y(t)$ for a local parameter $t$ such that $P = (x(0),y(0))$. To
explicitly compute with power series, we need to know where
($t$-adically) it is acceptable to truncate them.





\subsubsection{Precision: Cup product matrix.}
The first instance this problem arises is in the computation of the
cup product matrix. Since $v_t(\omega_j) = 2(g-j)-2,$ which is minimal
for $j = 2g-1$, we have $\min v_t(\omega_j) = -2g$. Thus it suffices
to compute each basis differential $\omega_k$ to a
precision of $t^{2g}$. Consequently, we compute $x(t)$ to a precision
of $t^{2(2g-1)}$ and $y(t)$ to a precision of $t^{2g-1}$.

\subsubsection{Precision: $\sum\Res(\alpha\int\beta)$.}
Now we consider $\sum_A\Res_A(\alpha\int\beta)$, where the sum is
taken over all points $A$ in the set $\{$poles of $\alpha$,
Weierstrass points of $C\}$. We begin by looking at the expansion for
$\alpha = \phi^*\omega - p\omega$.

If $A$ is a non-Weierstrass point defined over $\Q_p$, then computing
$\alpha$ in local coordinates is unnecessary, as $\alpha$ just has a
simple pole at $A$. So $\Res_A(\alpha\int\beta) = \Res_A \alpha$, and
we may simply read off the residue:\begin{itemize}\item $A=P$
  Teichm\"{u}ller: $1-p$ \item $A=P$ not Teichm\"{u}ller: $-p$
\item $A=-P$ Teichm\"{u}ller: $p-1$
\item $A=-P$ not Teichm\"{u}ller: $p$.\end{itemize}

In the case where $A$ is non-Weierstrass and defined over an extension
of $\Q_p$, we have to compute Coleman integrals and local coordinates,
so we must study the precision of both, as given by the following corollary of Proposition~\ref{prop:tiny precision}:
\begin{corollary}Let $A$ be non-Weierstrass, defined over degree $p$ extension $k'$ 
  of $\Q_p$. Let $U_A$ denote the residue disc of $A$, and let $B$ be a non-Weierstrass point in $U_A$ defined over $\Q_p$. Suppose a working precision of $n$ $p$-adic digits (so that $A$ has precision $pn$ in a uniformizer $\pi$).
  Let $\beta$ be written in terms of the local coordinate $(x(t),y(t))$ at $B$, so that $\beta = h(t)dt$ with $h(t)$ truncated modulo $t^{mp}$. Then the residue
  $\sum_{S \in U_A} \Res(\alpha \int \beta)$ has $\min\{n, pm + 1- \lfloor \log_p(pm+1)\rfloor\}$ digits of precision.
  
  \end{corollary}
\begin{proof}If $A$ is non-Weierstrass and defined over an extension $K$ of $\Q_p$, then the contribution of $\alpha$ to the residue calculation just depends on the disc of $A$:  $+1$ if $A$ is in $U_P$, and $-1$ if $A$ is in $U_{-P}$, with the constant of integration fixed in each residue disc. However, in this case, since we have multiple poles in each residue disc, we must compute $\Res(\alpha(\int_Z \beta)$, where the integral of $\beta$ is definite, taking into account a constant of integration chosen for each residue disc. More precisely: suppose we are working in the residue disc of $P$, and say $A$ is defined over a degree $p$ extension of $\Q_p$. Note that we must compute the local coordinates $(x(t),y(t))$ at $P$ with a precision of at least $t^{pm}$. As the interpolation from $P$ to $A$ is linear, we merely make a linear substitution \begin{align*}x(t) &:= x((x(A)-x(P))t)\\
    y(t) &:= y((x(A)-x(P))t).\end{align*} This new $x(t),y(t))$ is
  used to compute the tiny integral of $\beta$ from $P$ to $A$, the
  result of which has precision $\min\{n, pm + 1- \lfloor \log_p(pm+1)\rfloor\}$. Taking the trace from $K$ to $\Q_p$ accounts for
  the other poles of $\alpha$ in the disc of $P$.\end{proof}

Finally, in the case where $A$ is a finite Weierstrass point, we have to compute in the
local coordinates of $A$. (Note that we need not compute the residue at $(0,0)$ if on the curve or at infinity.)

\begin{proposition}Let $\alpha$ be above and let $A$ be a finite
  Weierstrass point not equal to $(0,0)$. Let $(x(t),y(t))$ represent
  the local coordinates at $A$. Then to compute
  $\Res(\alpha\int\beta)$ with $n$ digits of $p$-adic precision, we
  compute $(x(t),y(t))$ to $t^{2pn - p -3}$.\end{proposition}

\begin{proof}We have \begin{align*}\alpha &= \phi^*\omega - p\omega \\
    &= \frac{y(P)px^{p-1}{dx}}{\phi(y)(x^p-x(P))} - \frac{py(P)dx}{y(x-x(P))},\end{align*}
  where $$\frac{1}{\phi(y)}=y^{-p}\sum_{i=0}^{\infty}\binom{-1/2}{i}\frac{(f(x^p)-f(x)^p)^i}{f(x)^{pi}}.$$
  For $\Res(\alpha\int\beta)$ to have $n$ digits of $p$-adic
  precision, we must compute $n$ terms of the binomial expansion of
  $\frac{1}{\phi(y)}$.

  Recall that for a finite Weierstrass point $(a,0)$, we have \begin{align*}x(t) &= a + \frac{1}{g(a)}t^2 + O(t^4)\\
    y(t) &= t,\end{align*} where $g(x) = \frac{f(x)}{x-a}$. Note that
  by hypothesis, $a \neq 0$. We compute the $t$-adic
  valuation of $\alpha$:
  \begin{align*}v_t(\alpha) &= v_t(\phi^*\omega) \quad\text{since $\omega$ only contributes higher-order terms} \\
    &= v_t\left(\frac{y(P)px^{p-1}{dx}}{\phi(y)(x^p-x(P))}\right) \\
    &= 1 + v_t\left(\frac{1}{\phi(y)}\right) \quad(x^p \neq x(P))\\
    &= 1 - pv_t(y) + (n-1)v_t\left(\frac{f(x^p)-f(x)^p}{f(x)^p}\right)\\
    &= 1 - p + (n-1)\begin{cases} & (2-2p), \quad\text{if}\; v_t(f(x^p)-f(x)^p) > 0\\
      & -2p, \quad\quad\text{else}\end{cases}\end{align*}
  Thus we have  $$v_t\left(\alpha\right) = \begin{cases} & p - 2pn + 2n-1, \quad\text{if}\; v_t(f(x^p)-f(x)^p) > 0\\
    & p - 2pn+1, \quad\quad\text{else}.\end{cases}$$

  As $p -2pn + 2n-1 \geq p - 2pn +1 $ for $n \geq 1$, we have
  $v_t(\alpha) \geq p - 2pn + 1$. Set $m=2pn - p -1.$ Since we want
  $\Res(\alpha\int\beta)$, we need $v_t(\alpha\int\beta) \geq -1$, so
  we must compute $\beta$ to at least $t^{m-2}$. To get this
  precision, we must in turn compute with $x(t),y(t)$ to this
  precision.  \end{proof}

\section{Examples}\label{ex}
Here we provide some examples of our algorithms.

\subsection{Local heights: genus 2, general divisors}\label{subs:genex}
Let $C$ be the genus 2 hyperelliptic curve $$ y^2 = x^5 - 23x^3 + 18x^2 + 40x = (x - 4)(x - 2)x(x + 1)(x +
5)$$ over $\Q_{11}$, and let \begin{align*}D_1 &= (P) - (Q) \\
D_2 &= (R) - (S),\end{align*}
where $P = (-4,24), Q = (1,6), R = (5,30), S = (-2,12)$. We describe how to use Algorithm~\ref{generalalgorithm} to compute the local contribution at $p = 11$.

We see that \begin{align*}D_1^+ = \div\left(\frac{x-x(P)}{x-x(Q)}\right), &\quad D_1^-= [(P)-(-P)] + [(-Q)-(Q)], \\
D_2^+ = \div\left(\frac{x-x(R)}{x-x(S)}\right), &\quad D_2^- = [(R)-(-R)] + [(-S)-(S)].\end{align*}

Furthermore, we have \begin{align*}\frac{1}{4}h_{11}(D_1^+,D_2^+) &= \frac{1}{4}\log\left(\frac{x-x(P)}{x-x(Q)}(D_2^+)\right)\\
&= \frac{1}{2}\log\left(\left(\frac{x(R)-x(P)}{x(R)-x(Q)}\right)\left(\frac{x(S)-x(P)}{x(S)-x(Q)}\right)^{-1}\right)\\
&=2 \cdot 11 + 9 \cdot 11^{2} + 7 \cdot 11^{3} + 2 \cdot 11^{4} + O(11^{5}).\end{align*}

Now we compute, using Algorithm~\ref{antisymmetricalgorithm}, the contributions from antisymmetric heights (details of which are in Subsection~\ref{subs:specex}):
\begin{align*}h_{11}((P)-(-P),(-S)-(S)) &= 9 \cdot 11^{-1} + 5 + 6 \cdot 11 + 8 \cdot 11^{2} + 9 \cdot 11^{3} + 3 \cdot 11^{4} + O(11^{5})
 \\
h_{11}((P)-(-P),(R)-(-R)) &= 6 \cdot 11^{-1} + 10 + 7 \cdot 11 + 6 \cdot 11^{2} + 3 \cdot 11^{3} + 7 \cdot 11^{4} + O(11^{5})
\\
h_{11}((-Q)-(Q), (R)-(-R)) &= 8 \cdot 11^{-1} + 5 + 7 \cdot 11 + 10 \cdot 11^{2} + 3 \cdot 11^{3} + 7 \cdot 11^{4} + O(11^{5})
\\
h_{11}((-Q)-(-Q),(-S)-(S)) &= 11^{-1} + 8 + 7 \cdot 11 + 2 \cdot 11^{2} + 7 \cdot 11^{3} + 8 \cdot 11^{4} + O(11^{5}),\end{align*}
which gives
\begin{align*}\frac{1}{4}h_{11}(D_1^-,D_2^-) &= \frac{1}{4}(h_{11}((P)-(-P),(-S)-(S)) + \\
&\quad\;\quad h_{11}((P)-(-P),(R)-(-R)) + \\
&\quad\;\quad h_{11}((-Q)-(Q), (R)-(-R)) + \\
&\quad\;\quad h_{11}((-Q)-(-Q),(-S)-(S)))\\
&= 6 \cdot 11^{-1} + 7 + 4 \cdot 11 + 4 \cdot 11^{2} + 3 \cdot 11^{3} + 11^{4} + O(11^{5}).\end{align*}

Finally, we have \begin{align*}h_{11}(D_1,D_2) &= \frac{1}{4}h_{11}(D_1^+,D_2^+) + \frac{1}{4}h_{11}(D_1^-,D_2^-) \\
&= 6 \cdot 11^{-1} + 7 + 6 \cdot 11 + 2 \cdot 11^{2} + 4 \cdot 11^{4} + O(11^{5}).\end{align*}

\subsection{Local heights: genus 2, antisymmetric divisors}\label{subs:specex}
Keeping notation as in Subsection~\ref{subs:genex}, we describe in more detail how to use Algorithm~\ref{antisymmetricalgorithm} to compute the local contribution for one of the antisymmetric divisors: $$h_{11}((P)-(-P),(R)-(-R)).$$ For ease of notation, let us call these divisors $$D_P = (P)-(-P),\quad D_R = (R)-(-R).$$

With respect to the standard basis $\mathcal{B}$, the cup product matrix is
$$N = \left(\begin{tabular}{cccc}\label{cpm}
    0  &   0  &   0  &  $\frac{1}{3}$\\
    0  &   0  &   1  &   0\\
    0 &   $-1$  &   0 & $-\frac{23}{3}$\\
 $-\frac{1}{3}$ &    0  & $\frac{23}{3}$  &   0\end{tabular}\right).$$

Let $\nu_P$ be a differential with residue divisor $D_P$: we have $\nu_P = \frac{24dx}{y(x+4)}.$

We compute $\Psi(\nu_P)$ with respect to the basis $\{\omega_0,\omega_1, \Frob^n\omega_2, \Frob^n\omega_3\}$:
\begin{equation*}\Psi(\nu_P) = \left( \begin{array}{c}
8 \cdot 11^{-1} + 9 + 6 \cdot 11 + 3 \cdot 11^{2} + 7 \cdot 11^{3} + 11^{4} + O(11^5) \\
7 \cdot 11^{-1} + 1 + 4 \cdot 11^{2} + 2 \cdot 11^{3} + 8 \cdot 11^{4} + O(11^5) \\
7 + 9 \cdot 11 + 7 \cdot 11^{2} + 4 \cdot 11^{4} +8 \cdot 11^{5} + O(11^6)\\
2 + 2 \cdot 11 + 8 \cdot 11^{2} + 6 \cdot 11^{3} + 7 \cdot 11^{4} + 2 \cdot 11^{5} +O(11^6)
\end{array} \right).\end{equation*}

This implies that $$\int_{D_R}\eta_1 = 5 \cdot 11^{-1} + 6 + 3 \cdot 11 + 11^{3} + 7 \cdot 11^{4} + O(11^5).$$

To integrate $\nu_P$, we compute several quantities. Noting that $\alpha = \phi^*\nu_P-p\nu_P$ and that the $\Psi$ map is Frobenius equivariant, we have $$\Psi(\alpha) = \Psi(\phi^*\nu_P-p\nu_P) = \phi^*(\Psi(\nu_P)) -p\Psi(\nu_P).$$ In particular, this makes the computation of $\phi^*(\Psi(\nu_P))$ rather easy, as we have already computed $\Psi(\nu_P)$, and all that is left to do is multiply by the matrix of Frobenius.
We find that \begin{equation*}\Psi(\alpha) = \left( \begin{array}{c}
6 \cdot 11 + 5 \cdot 11^{2} + 2 \cdot 11^{4} + 9 \cdot 11^{5} + O(11^{6})\\
2 \cdot 11 + 10 \cdot 11^{2} + 8 \cdot 11^{3} + 6 \cdot 11^{4} + 2 \cdot 11^{5} + O(11^{6})\\
4 \cdot 11 + 6 \cdot 11^{2} + 2 \cdot 11^{3} + 11^{4} + 9 \cdot 11^{5} + O(11^{6})\\
3 \cdot 11 + 2 \cdot 11^{2} + 8 \cdot 11^{3} + 2 \cdot 11^{4} + 4 \cdot 11^{5} + O(11^{6})\end{array}\right).\end{equation*}

We wish $\beta$ to have residue divisor $D_R$, so let $\beta =\frac{30dx}{y(x-5)}.$ Then $$\Psi(\alpha) \cup \Psi(\beta) = 6 + 11^{2} + 9 \cdot 11^{4} + 5 \cdot 11^{5} + O(11^6).$$

To compute $\sum\Res(\alpha\int\beta)$, we must sum over all Weierstrass point and poles of $\alpha$. Recall that within a single residue disc, $\sum_A\Res_A(\alpha) =0$. Now computing the action of $\Psi$ on this differential is slightly more complicated, since instead of just two non-Weierstrass poles, we have $2p = 2\cdot 11$ non-Weierstrass poles: those points in the residue discs of $P$ and $-P$ with $x$-coordinate $\zeta_{11}^j(-4)^{1/11}$ (where $j = 0,\ldots, 10$).
This means we must work over the splitting field $L_{-4} = \Q_{11}(\zeta_{11}, (-4)^{1/11})$ of $x^{11}+4$ over $\Q_{11}$  to compute the local symbols. Since each set of $p$th roots is Galois conjugate, working over $L_{-4}$ yields $$\sum_j \langle \nu_1, \omega_i \rangle_{P_j} = \tr_{L_{-4}/\Q_{11}} (\langle \nu_1, \omega_i \rangle_{P_1}),$$ where $P_j$ is the point in the residue disc of $P$ with $x$-coordinate $\zeta_{11}^j (-4)^{1/11}$.
We have the following contribution from the disc of $P$:$$10 \cdot 11 + 9 \cdot 11^{2} + 4 \cdot 11^{3} + 3 \cdot 11^{4} + 4 \cdot 11^{5} + O(11^{6}),$$ and the total contribution from non-Weierstrass points is twice this, or
$$9 \cdot 11 + 8 \cdot 11^{2} + 9 \cdot 11^{3} + 6 \cdot 11^{4} + 8 \cdot 11^{5} + O(11^{6}).$$ Meanwhile, the sum of contributions from all Weierstrass discs is the following: $$11 + 4 \cdot 11^{3} + 6 \cdot 11^{4} + 11^{5} + O(11^{6}).$$

We compute the tiny integral $$\int_{R}^{\phi(R)}\nu_1= 8 \cdot 11 + 11^{2} + 8 \cdot 11^{3} + 2 \cdot 11^{5} + O(11^6).$$
Putting all of this together, we have $$h_{11}(D_P,D_R) = 6 \cdot 11^{-1} + 10 + 7 \cdot 11 + 6 \cdot 11^{2} + 3 \cdot 11^{3} + 7 \cdot 11^{4} + O(11^{5}).$$

As a consistency check, we also compute $h_{11}(D_R,D_P)$. Here we have $$\int_{D_P} \nu_R =2 + 11^{3} + 10 \cdot 11^{4} + 4 \cdot 11^{5} + O(11^{6})$$
and
$$\int_{D_P} \eta_R = 5 \cdot 11^{-1} + 2 + 3 \cdot 11 + 4 \cdot 11^{2} + 8 \cdot 11^{3} + 2 \cdot 11^{4} + O(11^{5}),$$
which gives
$$h_{11}(D_R,D_P) = 6 \cdot 11^{-1} + 10 + 7 \cdot 11 + 6 \cdot 11^{2} + 3 \cdot 11^{3} + 7 \cdot 11^{4} + O(11^{5}),$$ illustrating symmetry of the local height pairing.


\subsection{Global heights: genus 1}

We give an example of our implementation in genus 1, which allows for comparison of global heights via the algorithm of Mazur-Stein-Tate.

Let $C$ be the elliptic curve $$y^2 = x^3 - 5x,$$ with $Q = (-1,2),R = (5,10)$, so that $$(Q)-(-Q) = (R)-(-R) = \left(\frac{9}{4},-\frac{3}{8}\right) =:P.$$

We compute the $13$-adic height of $P$:
\begin{itemize} \item Above 13, the local height $h_{13}((Q)-(-Q),(R)-(-R))$ is $$2\cdot 13 + 6\cdot 13^2 + 13^3 + 5\cdot 13^4 + O(13^5).$$
\item Away from 13, the only nontrivial contribution is at 3, which is $2\log 3$ (by work of M\"{u}ller).
\item So the global 13-adic height is $12\cdot 13 + 4 \cdot 13^2 + 10 \cdot 13^3 + 9 \cdot 13^4 + O(13^5)$.
\end{itemize}
We compare this to Harvey's implementation \cite{harvey:heights} of the Mazur-Stein-Tate algorithm in Sage:
\begin{verbatim}
sage: C = EllipticCurve([-5,0])
sage: f = C.padic_height(13)
sage: f(C(9/4,-3/8)) + O(13^5)
12*13 + 4*13^2 + 10*13^3 + 9*13^4 + O(13^5)
\end{verbatim}

\section{Future work}\label{future}
Below we discuss some natural questions arising from our work.

\subsection{Global height pairings}
Ultimately, we would like to compute the global height pairing. To do
so, we would again require $C$ to be a curve over a number field $K$
with good reduction at each place $v$ dividing $p$. We would also need
a continuous id\`{e}le class character $\chi: \A_K^*/K^* \lra \Q_p$
and a splitting $H^1(C/K_v) = H^{1,0}(C/K_v) \oplus W_v$ for each
place $v$ dividing $p$. Computing the local heights at those primes
$v$ away from $p$ and those above $p$, the global height would then be
the sum of all local heights. When $K = \Q$, the recent Ph.D. thesis of M\"{u}ller \cite{muller:thesis} addresses these local heights away from $p$, and putting together our results, we are able to compute global heights, as shown in Section~\ref{ex}. It would be quite interesting to extend these computations to number fields.

\subsection{Optimizations}
In another direction, it is also of interest to optimize the present algorithm. Currently,
the most expensive step is in computing the Laurent series expansion
of $\alpha$ in the various Weierstrass local coordinates to reasonably
high precision. As we are just interested in the residue of
$\alpha\int\beta$, is there a way to make this faster?

\subsection{Comparison with the work of Mazur-Stein-Tate}
When the curve is elliptic, we are able to compare our algorithm for the height pairing with the algorithm
of~\cite{mazur-stein-tate}, as in Section~\ref{ex}. We note that we compute the height pairing
for divisors with disjoint support. It is obviously possible to
compute without this assumption by replacing one divisor by a linearly
equivalent one with this property. But it is also possible to extend
the method described in~\cite[\S 5]{Gross86}. This extended method can be
compared directly with the method of~\cite{mazur-stein-tate}, as the
height is just the height pairing of a divisor with itself. 